\documentclass{amsart}
\usepackage{amsmath}
\usepackage{graphicx} 
\usepackage{amsfonts} 
\usepackage{amssymb}
\usepackage{color}
\usepackage{mathrsfs}
\usepackage{epsfig}
\usepackage{subfig}
\usepackage{url}
\theoremstyle{plain}
\newtheorem{theorem}{Theorem}[section]

\newtheorem{lemma}{Lemma}[section]
\newtheorem{proposition}{Proposition}[section]

\theoremstyle{definition}
\theoremstyle{remark}
\numberwithin{equation}{section}

\newcommand{\R}{\mathbb{R}}

\newcommand{\Z}{\mathbb{Z}}

\newcommand{\calE}{{\mathcal E}}

\newcommand{\calL}{{\mathcal L}}
\newcommand{\calC}{{\mathcal C}}

\newcommand{\VHR}{{V}_{HR}}

\definecolor{viole}{RGB}{100,0,120}
\definecolor{darkor}{RGB}{220,110,0}

\usepackage{latexsym}

\title[Crystallization and discrete Gauss-Bonnet]{Crystallization in two dimensions and a discrete Gauss-Bonnet theorem}

\author[L. De Luca]
{L. De Luca}
\address[Lucia De Luca]{ Zentrum Mathematik - M7, Technische Universit\"at  M\"unchen, Boltzmannstrasse 3, 85748 Garching, Germany}
\email[L. De Luca]{deluca@ma.tum.de}

\author[G. Friesecke]
{G. Friesecke}
\address[Gero Friesecke]{ Zentrum Mathematik - M7, Technische Universit\"at  M\"unchen, Boltzmannstrasse 3, 85748 Garching, Germany}
\email[G. Friesecke]{gf@ma.tum.de}

\begin{document}

\begin{abstract}
We show that the emerging field of discrete differential geometry can be usefully brought to bear on crystallization problems. In particular, we give a simplified proof of the Heitmann-Radin crystallization theorem \cite{HR}, which concerns a system of $N$ identical atoms in two dimensions interacting via the idealized pair potential $V(r)=+\infty$ if $r<1$, $-1$ if $r=1$, $0$ if $r>1$. This is done by endowing the bond graph of a general particle configuration with a suitable notion of {\it discrete curvature}, and appealing to a {\it discrete Gauss-Bonnet theorem} \cite{Knill1} which, as its continuous cousins, relates the sum/integral of the curvature to topological invariants. This leads to an exact geometric decomposition of the Heitmann-Radin energy into (i) a combinatorial bulk term, (ii) a combinatorial perimeter, (iii) a multiple of the Euler characteristic, and (iv) a natural topological energy contribution due to defects. An analogous exact geometric decomposition is also established for soft potentials such as the Lennard-Jones potential $V(r)=r^{-6}-2r^{-12}$, where two additional contributions arise, (v) elastic energy and (vi) energy due to non-bonded interactions. 
\end{abstract}

\maketitle



\section{Introduction} 
At low temperature, atoms and molecules typically self-assemble into crystalline order. The challenge to derive this fundamental phenomenon from a mathematical model of the interatomic interactions is known as the {\it crystallization problem}. In the limit of zero temperature and long time, observed configurations are expected to correspond to minimizers of the interaction energy, and so the crystallization problem amounts to proving that energy minimizers exhibit crystalline order. 


Our goal in this paper is to understand an important initial result on crystallization in two dimensions in a new way, by introducing and exploiting a {\it discrete differential geometry} viewpoint. We hope that this approach will also aid future progress on the many open crystallization problems, including 3D ones. In a companion paper \cite{DF2}, our approach will be used for developing a rigorous understanding of basic defects occuring at non-minimal but low energy for the models studied here. 

\begin{figure}[http!]
\begin{center}
\includegraphics[width=0.25\textwidth]{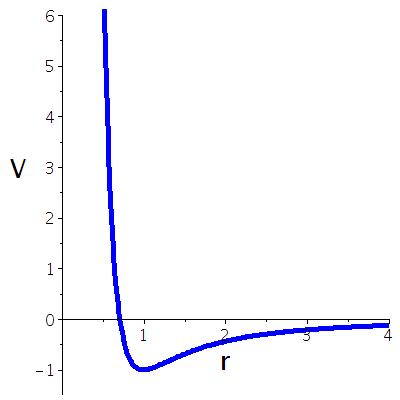} \hspace*{5mm} 
\includegraphics[width=0.25\textwidth]{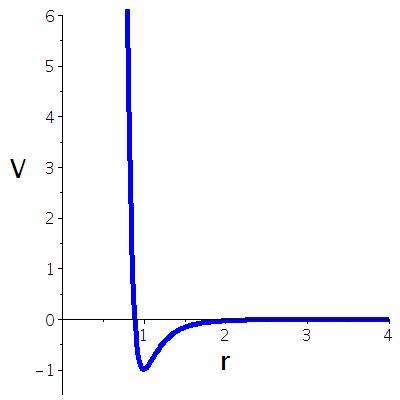} \hspace*{5mm}
\includegraphics[width=0.25\textwidth]{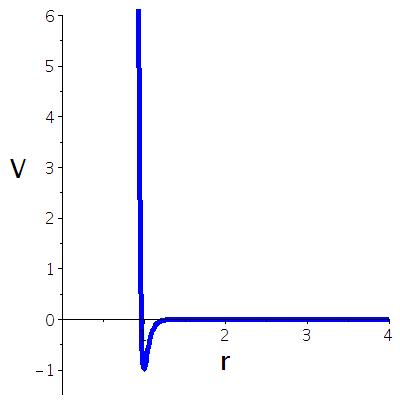}
\end{center}
\caption{The Lennard-Jones potential \eqref{LJ} for $p=2$ (left), $6$ (middle) and $20$ (right).} 
\label{F:potentials}
\end{figure}

The important crystallization result we are concerned with is due to Heitmann and Radin \cite{HR}. One starts from the prototypical Lennard-Jones energy for a system of $N$ identical atoms with positions $x_1,..,x_N\in\R^d$, 
\begin{equation}\label{energy}
  \calE_{_{V}}(x_1,..,x_N) = \sum_{i<j} V(|x_i-x_j|)
\end{equation}
where $V$ is the Lennard-Jones $(p,2p)$ potential 
\begin{equation} \label{LJ} 
     V_{p,2p}(r) = r^{-2p} - 2 r^{-p} \;\; (p>0).
\end{equation}
(see Figure \ref{F:potentials}). Here we have normalized the optimal interparticle distance to $1$ and the associated potential energy to $-1$.
This energy numerically exhibits crystallization. See \cite{WalesDoye} for a list of numerical minimizers in three dimensions for $p=6$ and $N$ up to $110$. The model \eqref{energy}--\eqref{LJ} plays an important role in the physics literature on large systems of atoms above and below the crystallization temperature (see e.g. \cite{AllenTildesley, TroianEtAl, WalesDoye, CameronVandenEijnden}). For further information about model energies for many-atom systems, including highly accurate quantum-mechanical ones, we refer to the recent survey article \cite{FrieseckeTheil}.
As regards mathematical results, even in dimension $d=2$ crystallization of minimizers has not been proved rigorously (see the end of this introduction for further discussion of this point), but Heitmann and Radin established the beautiful result that when the potential is simplified to
\begin{equation} \label{HR}
  \VHR(r) = \begin{cases} + \infty, & r<1 \\  - 1, & r=1 \\ 0, & r>1, \end{cases}
\end{equation}
then all minimizers are up to rotation and translation subsets of the triangular lattice
\begin{equation} \label{lattice}
   \calL = \{i {\bf e} + j {\bf f} \, : \, i,j\in\Z\}, \;\; {\bf e}=\begin{pmatrix} 1 \\ 0 \end{pmatrix}, \;\; {\bf f} = \begin{pmatrix} 1/2 \\ \sqrt{3}/2 \end{pmatrix}.
\end{equation}
We note that the Heitmann-Radin potential is the large-$p$ limit of the Lennard-Jones potential, i.e. $\lim_{p\to\infty} V_{p,2p}(r)=V_{HR}(r)$; physically, one is taking a ``brittle limit'' in which the width of the well of the pair potential is compressed to zero and bonds immediately break upon increasing the interparticle distance.
For finer information about the minimizers of \eqref{energy}, \eqref{HR} such as their asymptotic shape see \cite{AFS, Schmidt, DPS}, and for an analogous result in the presence of a three-body potential favouring bond angles of $2\pi/3$ see \cite{MaininiStefanelli}.

Roughly, our new approach to understanding the Heitmann-Radin result is as follows. 
First, endow the bond graph of a particle configuration with a suitable notion of discrete (combinatorial) curvature which vanishes on the expected ground state graphs. Second, ``integrate'' over the whole configuration and use a discrete version of the Gauss-Bonnet theorem from differential geometry to obtain a purely topological contribution to the total energy. Third, show that this contribution as well as the bulk energy (of which more below) can be eliminated, and minimize just the remaining energy contributions. This strategy can be summarized in the following diagram:
$$
   \mbox{particle configuration} \longrightarrow \mbox{bond graph} 
   \longrightarrow \mbox{curvature} \longrightarrow \mbox{simplified energy}. 
$$
Finding a suitable notion of curvature which localizes a relevant part of the energy is a long story, told in Section \ref{S:DiffGeo}. As turns out, we will work with the
combinatorial Puiseux curvature (see \cite{Knill1} or \eqref{Puiseux} below) of the triangulated bond graph, i.e. the graph obtained by adding extra edges until all non-triangular faces have been decomposed into triangles. A recent discrete Gauss-Bonnet theorem by Knill \cite{Knill1} says that this curvature integrates to a multiple of the Euler characteristic of the bond graph (see also Theorem \ref{T:GB} below, which extends this result to irregular boundaries as needed here). This leads to the following exact expression for the energy \eqref{energy}, \eqref{HR} of an arbitrary configuration $X=(x_1,..,x_N)$ (see Theorem \ref{T:geo} below): 
\begin{equation} \label{intro:decomp1}
    \calE_{_{V_{HR}}}(X) = -3N + {P}(X) + \mu(X) + 3\chi(X),
\end{equation}
where $\chi(X)$ is the Euler characteristic of the bond graph and $\mu(X)$ is a natural {\it defect measure}, namely the number of additional edges due to triangulation
(see \eqref{mu}). For nice configurations, $\mu$ can be expressed in terms of the original bond graph as 
$$
    \mu(X)= \sharp \, \mbox{quadrilaterals} + 2 \, \sharp\, \mbox{pentagons} + 3 \, \sharp\,  \mbox{hexagons} + ...
$$  
The remaining energy contribution ${P}(X)$ in \eqref{intro:decomp1} is a {\it combinatorial perimeter}, defined on nice configurations as the number of boundary edges and extended to irregular configurations as suggested by geometric measure theory \cite{Federer}, namely by counting ``wire edges'' twice (see Figure \ref{F:edge} and \eqref{peri}). 

The Heitmann-Radin proof of crystallization relies on the following
remarkable but somewhat mysterious inequality due to Harborth \cite{Harborth}: for configurations with interparticle distance $\ge 1$ and simply closed polygonal boundary,
$$
     \calE_{_{V_{HR}}}(X) \ge \calE_{_{V_{HR}}}(X\backslash\partial X) - 3 \, \sharp \, \partial X + 6.
$$
The proof uses an ``elementary'' lower bound on the inner angle $\alpha(x)$ between incoming and outgoing boundary edge at the boundary vertex $x$ (see Figure \ref{F:Puiseux}),
\begin{equation}\label{harborth2}
     \alpha(x) \ge \Bigl(\sharp(\mbox{interior edges emanating from }x) + 1\Bigr) \frac{\pi}{3}.
\end{equation} 
Our simplified form of the Harborth inequality is 
\begin{equation} \label{FDeLuca2}
     (P+\mu)(X) \ge (P+\mu)(X\backslash\partial X) + 6
\end{equation}
(see Lemma \ref{L:lower}). The underlying ``elementary'' bound \eqref{harborth2} will be seen to have the following differential-geometric meaning: {\it the combinatorial Puiseux curvature as introduced by Knill \cite{Knill1} of the boundary of the bond graph is pointwise bounded from below by the euclidean Puiseux curvature of the associated polygon in the plane}. See Proposition \ref{P:curvaturebound}.
The global inequality \eqref{FDeLuca2} will be derived from this curvature bound by integration over the boundary and using Gauss-Bonnet. 

Another interesting outcome of the differential-geometric viewpoint is a generalization of the decomposition \eqref{intro:decomp1} to soft potentials such as \eqref{LJ}. In this case, edges of the bond graph are taken to correspond to particle pairs whose distance lies in a suitable neighbourhood $[\alpha,\beta]$ of the optimal distance $r=1$, and the energy \eqref{energy} is shown to satisfy
\begin{equation} \label{intro:decomp2}
   \calE_V(X) = -3N 
+ P(X) + 3\chi(X) + \mu(X) + \calE_{e\ell}(X) + \calE_{nbond}(X),
\end{equation}
where $\calE_{e\ell}$ is the elastic energy of the bonds which vanishes for optimal bondlengths $r=1$ (see \eqref{Eel} in Section \ref{S:energy}), and $\calE_{nbond}$ the energy due to non-bonded interactions (see \eqref{Enbond}). We emphasize that \eqref{intro:decomp2} holds for {\it arbitrary} configurations. Hence it is potentially useful not just for the study of minimizers, but also for studying crystallization at finite temperature or analyzing the molecular dynamics of crystal formation.

What can be said about minimizers in the soft potential case? One might expect that when $V$ is ``close'' to \eqref{HR}, minimizing configurations are ``close'' to subsets of a suitable lattice, up to a boundary layer of $o(N)$ particles. The deep insight of Theil \cite{Theil} (see \cite{ELi} for an extension and \cite{FlatleyTheil, FTTT} for a computer-aided generalization to a class of three-body potentials in 3D) that $N$-particle lattice subsets achieve the optimal asymptotic energy per particle in the limit $N\to\infty$ is very suggestive of such a behaviour of the positions. Unfortunately, the methods introduced here are insufficient; in particular, the pointwise bound of discrete Puiseux curvature by euclidean Puiseux curvature breaks down, and hence so does our differential-geometric proof of \eqref{FDeLuca2}.  
\section{Discrete differential geometry of particle configurations} \label{S:DiffGeo}
In this section we are concerned with {\it arbitrary} $N$-particle configurations in the plane. We endow them with additional mathematical structure (graph; metric; discrete curvature) and derive a discrete Gauss-Bonnet theorem by adapting recent work of Knill \cite{Knill1} to general triangular graphs.
\subsection{Configurations}
By a {\it configuration} of a system of $N$ identical particles in two dimensions we mean a set $X = \{x_1,..,x_N\}\subset\R^2$ of mutually distinct particle positions $x_i\in\R^2$.
\subsection{Bond graph} 
\begin{figure}[http!]
\begin{center}
\includegraphics[width=0.3\textwidth]{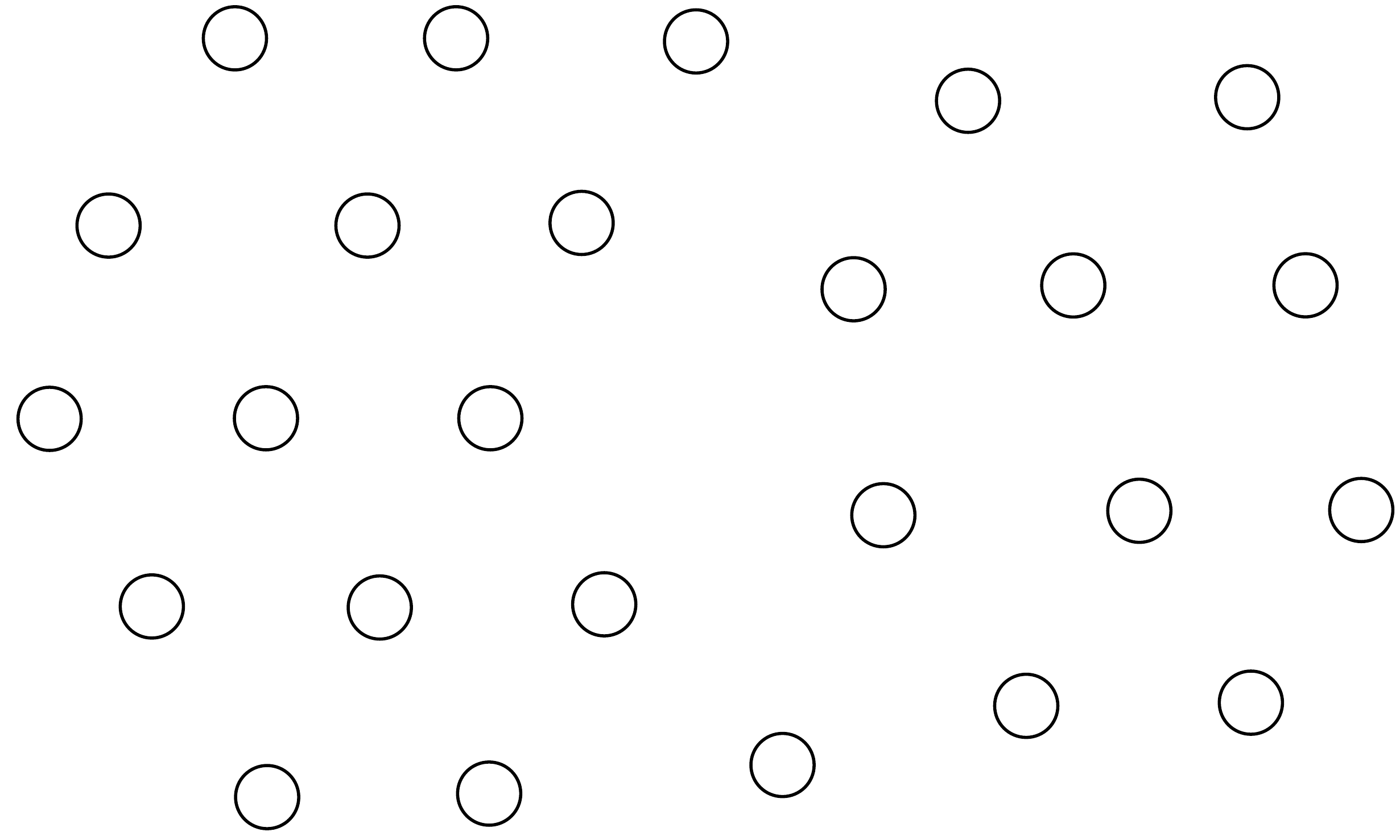} \hspace*{1cm} 
\includegraphics[width=0.3\textwidth]{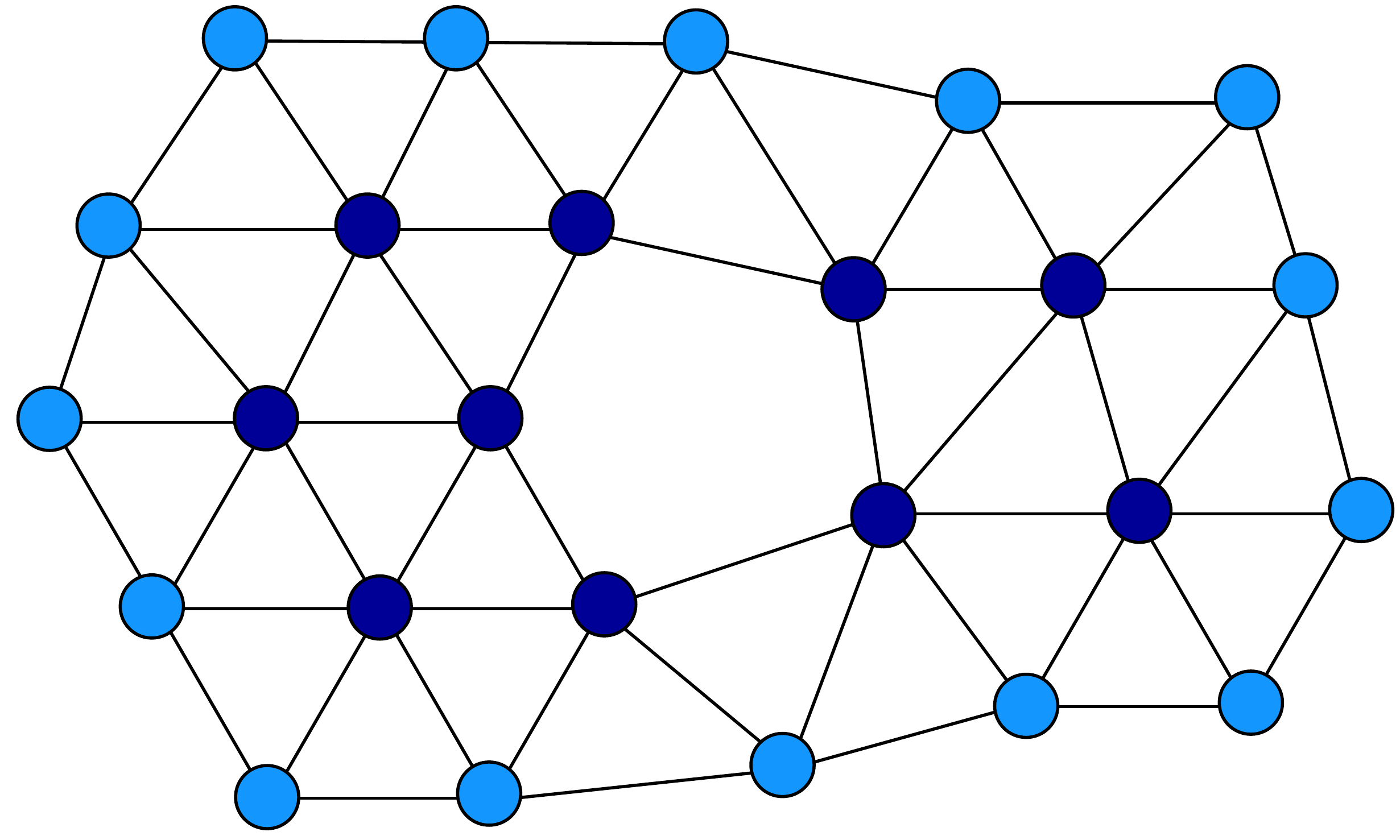}
\end{center}
\caption{A particle configuration and its bond graph. 
The colouring indicates boundary particles.} 
\label{F:bond}
\end{figure}

Discard now any detailed features of the pair potential $V$ and just assume that it has a unique minimum, say at $r=1$. We say that two particles $x$ and $y$ in $X$ are linked by an {\it edge}, or {\it bond}, if their mutual distance lies in a suitable neighbourhood, or {\it bond range}, $[\alpha,\beta]$ around the optimal distance $r=1$. Fix such a bond range $[\alpha,\beta]$, with $0<\alpha\le 1\le \beta$; the associated set of edges is  
\begin{equation}\label{edge}
  E := \{ \{x,y \} \, : \, x,y\in X, \, \alpha\le |x-y|\le \beta \}.
\end{equation}
We call the graph $(X,E)$ the {\it bond graph} of the configuration. See Figure \ref{F:bond}.

Simple geometric considerations show that the bond graph is a {\it planar graph} (i.e., for any two different edges $\{x,y\}$ and $\{x',y'\}$, the corresponding line segments $[x,y]$ and $[x',y']$ do not cross) provided 
\begin{equation}\label{planar1} 
             \beta < \sqrt{2}\, d_{min},
\end{equation}
where $d_{min}$ is the minimum interparticle distance $\min \{ |x-y| \, : \, x,y\in X, \, x\neq y\}$.  

Condition \eqref{planar1} can be fulfilled by a suitable choice of the maximum bondlength $\beta$ if and only if the particle configuration satisfies the mild minimal distance condition 
\begin{equation}\label{planar2}
             d_{min} > \frac{1}{\sqrt{2}}.
\end{equation}
We remark that if a particle pair achieved equality in \eqref{planar2}, its
standard Lennard-Jones (6,12) energy would be +48, a value almost two orders of magnitude above the binding energy $V(\infty)-V(1)=1$. Hence violations of \eqref{planar2} should never occur for planar ground state configurations of physical systems; for a rigorous proof in case of the Lennard-Jones (6,12) energy see \cite{DF2}. In the sequel we assume that \eqref{planar1}--\eqref{planar2} are satisfied.    
\subsection{Interior and boundary} \label{S:IandB}
\begin{figure}[http!]
\begin{center}
\includegraphics[width=0.5\textwidth]{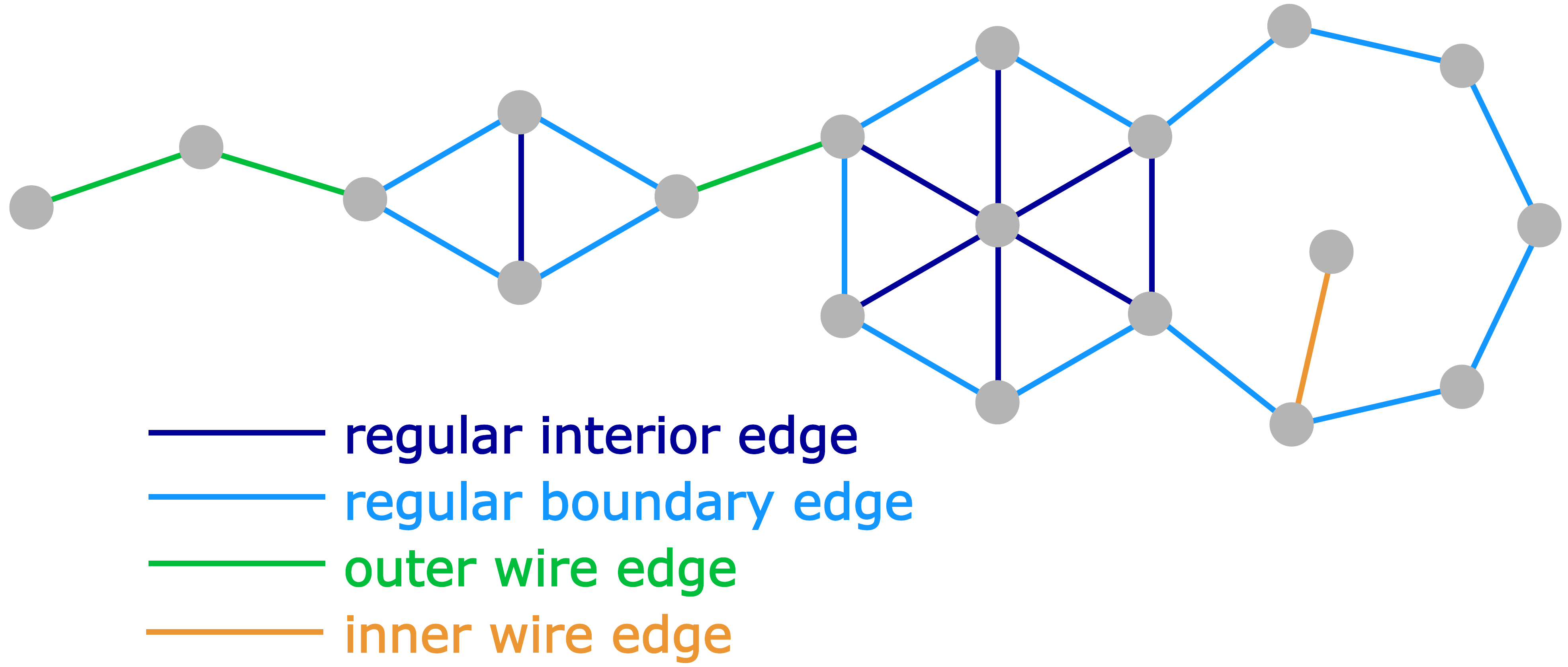}
\end{center}
\caption{Classification of edges} \label{F:edge}
\end{figure}

In order to study crystallization, it will be very useful to distinguish between ``interior'' and ``boundary'' particles and edges. To this end we first introduce {\it faces.} 
By a face $f$ we mean any open and bounded subset of $\R^2$ which is nonempty, does not contain any point $x\in X$, and whose boundary is given by a cycle, i.e., $\partial f = \cup_{i=1}^k [x_{i-1},x_i]$ for some points $x_0,x_1,..,x_k=x_0\in X$ with $\{x_{i-1},x_i\}\in E$.\footnote{The points $x_0,..,x_{k-1}$ do not need to be pairwise distinct, as a face might contain ``inner wire edges'' (see Figure \ref{F:edge}). Note also that for non-connected graphs, our definition differs slightly from standard conventions because ring-shaped regions bounded by two cycles are not faces. This has the advantage that the Euler characteristic stays unchanged under triangulation.} 
\\[2mm]
Edges can now be classified into four types (see Figure \ref{F:edge}). We say that an edge is 
\begin{itemize}
\item a {\it regular interior edge} if it lies on the boundary of two faces, 
\item a {\it regular boundary edge} if it lies on the boundary of precisely one face and that of the complement of the union of all faces,
\item an {\it outer wire edge} if it does not lie on the boundary of any face, 
\item an {\it inner wire edge} if it lies on the boundary of precisely one face but not on that of the complement of the union of all faces. 
\end{itemize}

The regular interior and inner wire edges are called {\it interior edges}, and the regular boundary and outer wire edges are called {\it boundary edges}. Finally we are in a position to distinguish mathematically between interior and boundary {\it particles} (see Figure \ref{F:bond}). We say that a particle is 
\begin{itemize}
\item an {\it interior particle} if it is incident only with interior edges, and
\item a {\it boundary particle} if it is isolated or incident with a boundary edge. 
\end{itemize}
The set of boundary particles will in the sequel be denoted by $\partial X$, and the set of interior particles by $\mbox{int}\, X$. 
\subsection{Perimeter}
As explained in the Introduction, an important role in our analysis is played by the {\it combinatorial perimeter} of the bond graph. For general configurations we define it as follows:
\begin{equation} \label{peri}
   P(X) := \sharp \mbox{(regular boundary edges in $E$)} 
            + 2 \sharp \mbox{(outer wire edges in $E$)}.
\end{equation}
We note that $P$ is additive over connected components, and coincides for connected graphs with edge length 1 with the geometric perimeter of the set $\tilde{X}$ in the plane obtained by taking the union of vertices, edges, and faces of the bond graph $(X,E)$. Note that the geometric perimeter, defined as the infimum of the length of simply closed smoooth curves whose interior contains $\tilde{X}$, naturally counts wire edges twice. 
\subsection{Which curvature relates to energy minimisation?}
As mentioned in the Introduction, it is not obvious which notion of curvature will be the most fruitful for our purposes. Here are some desiderata. 
\begin{enumerate}
\item {\it Curvature should be well-defined for general (irregular) configurations.} 
\item {\it For lattices, zero curvature should single out the energy-minimizing lattices.} 
\item {\it Non-topological defects such as elastic deformation, vacancies, or flat boundaries should not contribute to total curvature.}
\end{enumerate}
Condition (1) will lead to considering, along with the bond graph, its triangulation, of which more later. Condition (3) suggests to work with a purely combinatorial notion. Condition (2) rules out ``universal'' notions such as {\it Gromov curvature} \cite{Gromov, Higuchi}
\begin{equation}\label{Gromov}
    K_{Gr}(x) = 1 - \frac12 \mbox{(edge degree of $x$)} + \sum_{y\in S_1(x)} \frac{1}{\mbox{(face degree of $y$)}},
\end{equation}
where here and below $S_1(x)$ denotes the unit sphere with respect to the graph metric around a point $x$, that is to say all $y\in X$ linked to $x$ by an edge. This is because Gromov curvature vanishes on {\it all} of the standard lattices such as the triangular, square, and hexagonal lattice.

Euler curvature 
\begin{equation}\label{Euler}
    K_{Eu}(x) = 1 - \frac12 \mbox{(edge degree of $x$)} + \frac13 \mbox{(face degree of $x$)}
\end{equation}
distinguishes between these lattices, and hence satisfies condition (2), but it violates condition (3), because it yields a positive contribution from each boundary atom. 

\subsection{Curvature}
What turns out to work is a discrete version of {\it Puiseux curvature}, applied not to the bond graph but the {\it triangulated bond graph}. Discrete Puiseux curvature was recently introduced and studied in the context of subdomains of the triangular lattice by Knill \cite{Knill1}, and can be thought of as a boundary-corrected version of Euler curvature, corrected in such a way that flat boundaries have curvature zero. 

\begin{figure}[http!]
\begin{center}
\includegraphics[width=0.3\textwidth]{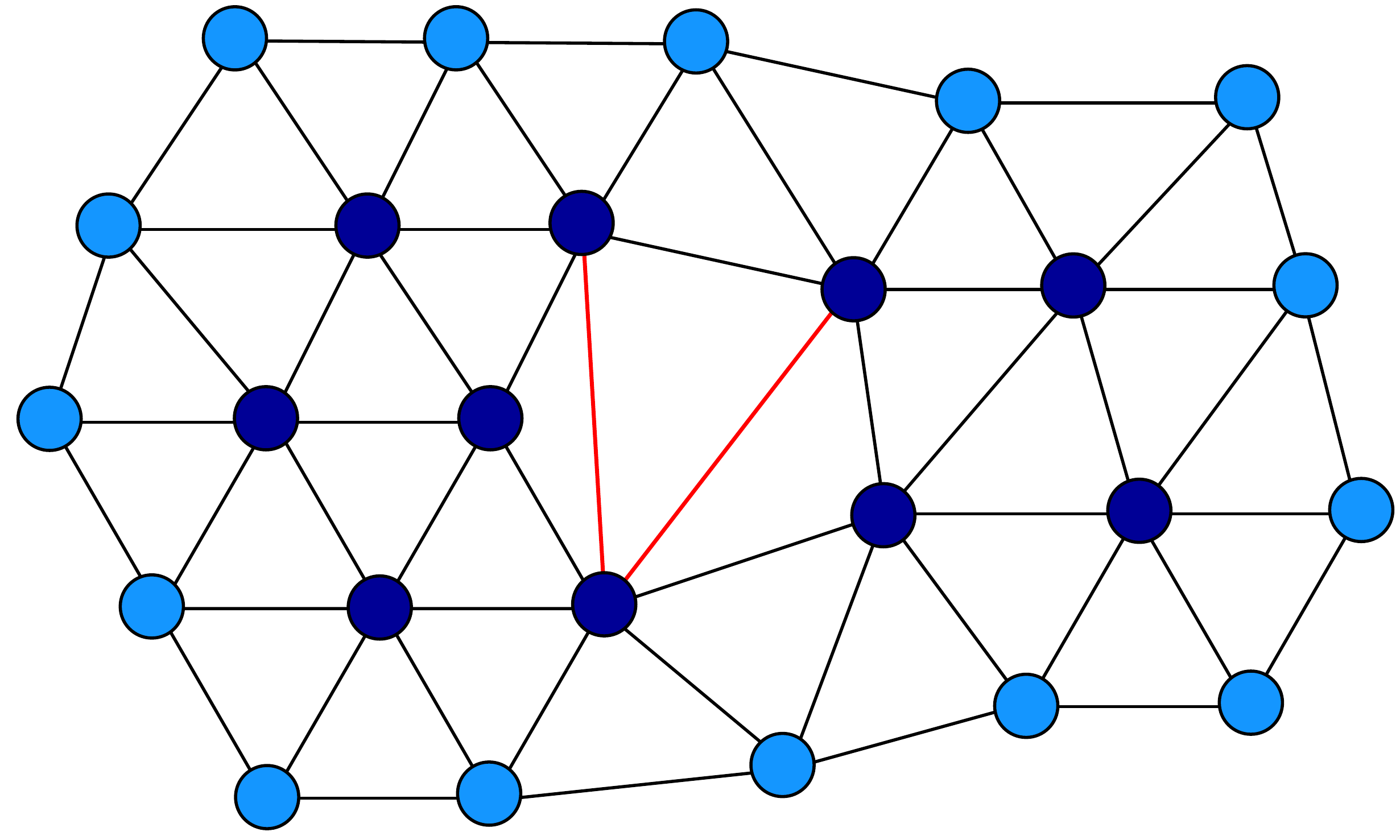}
\end{center}
\caption{The triangulated bond graph of the particle configuration from Figure \ref{F:bond}. Additional edges due to triangulation of faces with more than 3 edges are shown in red. 
} 
\label{F:triang}
\end{figure}

First we define the triangulated bond graph. For any face of the bond graph $(X,E)$ with more than 3 edges, say $k$ edges (where inner wire edges are counted twice), we add $k-3$ edges connecting not already connected vertices and not crossing each other. This yields a new graph $(X,\bar{E})$ all of whose faces are triangles and which we call the {\it triangulated bond graph}. See Figure \ref{F:triang}. 

\begin{figure}[http!]
\begin{center}
\includegraphics[width=0.95\textwidth]{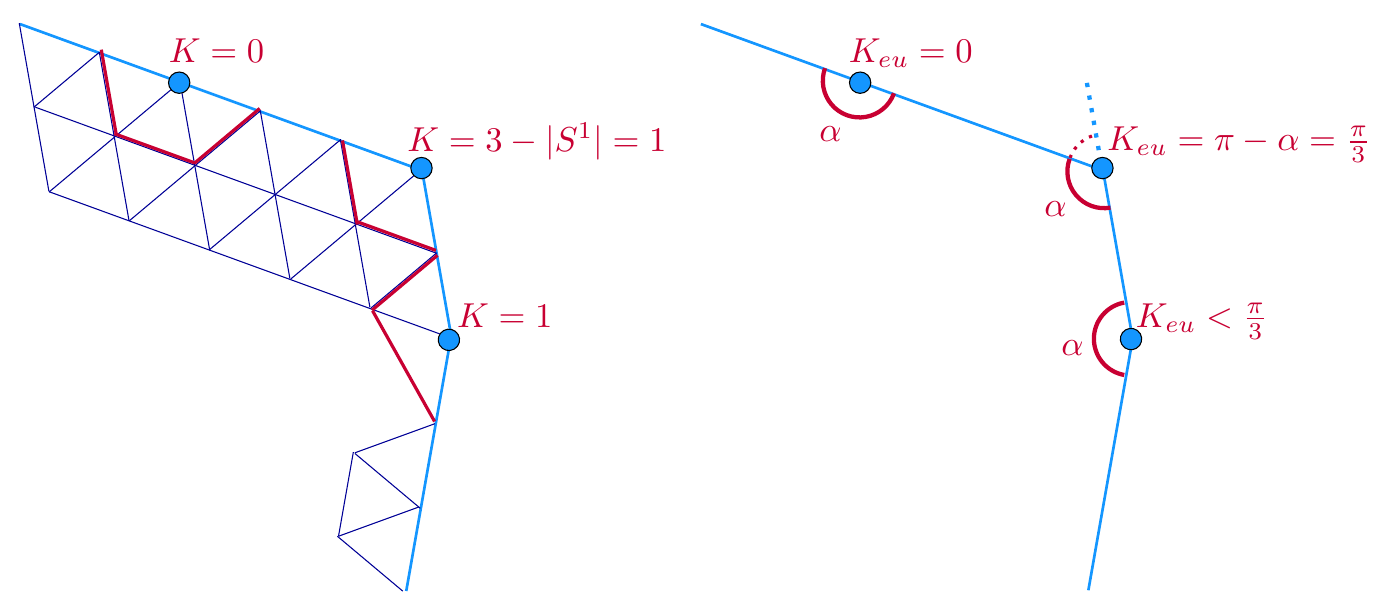}
\end{center}
\caption{Left: The combinatorial Puiseux curvature of a particle configuration at a boundary point is determined by the length of the unit sphere (drawn in red) in the graph metric. Right: The classical Puiseux curvature of the piecewise linear curve in the plane formed by the boundary edges is determined by the length of spheres (drawn in red) in the euclidean metric.} 
\label{F:Puiseux}
\end{figure}

We now define the {\it combinatorial Puiseux curvature} of a particle configuration $X$ at a point $x\in X$ (see \cite{Knill1}). Denote by $S_{1,{\calL}}$ the unit sphere around a point in the triangular lattice \eqref{lattice}. Let 
\begin{equation}\label{Puiseux}
    K(x) = \begin{cases} |S_{1,\calL}| - |S_1(x)|  \; = \; 6-|S_1(x)| & \mbox{if $x$ is an interior particle} \\
                        \mbox{$\frac12$}|S_{1,\calL}| - |S_1(x) \; = \; 3 - |S_1(x)| & \mbox{if $x$ is a boundary particle,}
           \end{cases}
\end{equation}
where ${S_1}(x)$ denotes the unit sphere with respect to the graph metric in the triangulated bond graph $(X,\bar{E})$ and $|S_1(x)|$ denotes the number of its edges. See Figure \ref{F:Puiseux}. 
To explain the meaning of \eqref{Puiseux} some remarks are in order. 

First, $\calL$ is the lattice expected to minimize the energy \eqref{energy}--\eqref{LJ}, and hence $K$ has the desired properties (2) and (3).
If interactions other than \eqref{energy}--\eqref{LJ} were under consideration, e.g. angular terms promoting the formation of graphene sheets and nanotubes as recently studied mathematically in \cite{MaininiStefanelli}, other reference lattices need to be used in \eqref{Puiseux}. 

Second, for interior particles $K(x)$ agrees (up to an overall normalization factor which is a matter of convention) with the Euler curvature \eqref{Euler}, because in a planar triangular graph both the edge degree and the face degree of an interior point $x$ agree with the number of edges in the unit sphere around $x$.

Third, for boundary particles, definition \eqref{Puiseux} is a graph-theoretic analogon of the classical euclidean Puiseux curvature of the boundary $\partial\Omega$ of a polygon $\Omega$ in $\R^2$, 
\begin{equation} \label{Keu}
    K_{eu}(x) = \underbrace{\mbox{$\frac12$}|S_1|_{eu}}_{=\pi} - \lim_{r\to 0} \frac{\mbox{length}(S_r(x)\cap\Omega)}{r} \;\; (x\in\partial\Omega),
\end{equation}
where $S_r(x)$ denotes the euclidean unit sphere of radius $r$ around $x$ and $|S_r(x)|_{eu}$ denotes its euclidean length. See Figure \ref{F:Puiseux}. The constant term in \eqref{Keu} for boundary particles is the euclidean length of a half-sphere in flat 2D euclidean space, just as the constant term in the second line of \eqref{Puiseux} is the length (w.r. to the graph metric) of the half-sphere in the triangular lattice. The second term in \eqref{Keu} just gives the interior angle. A deeper connection between combinatorial and euclidean Puiseux curvature will be given in Proposition \ref{P:curvaturebound}. 
\subsection{Discrete Gauss-Bonnet} \label{S:DGB}
The curvature $K(x)$ depends on the chosen triangulation, but -- as we shall now see -- its sum over all $x\in X$ does not; it is a topological invariant. The usefulness of $K(x)$ for crystallization problems stems from this property, as well as from the fact that it is naturally related to the energy \eqref{energy} (see \eqref{decomp1} below).  
Recall that the Euler characteristic of a planar graph with $v_0$ vertices, $v_1$ edges and $v_2$ faces is defined as 
\begin{equation} \label{Euchar}
   \chi(X) := \sum_{k=0}^2 (-1)^k v_k.
\end{equation}
\begin{theorem} \label{T:GB} (Discrete Gauss-Bonnet) Let $(X,\bar{E})$ be any planar triangular graph. Then 
\begin{equation} \label{GBrough}
   \sum_{x\in X} K(x) = 6 \chi(X) + 3\Bigl( P(X) - \sharp \partial X\Bigr).
\end{equation} 
In particular, when $(X,\bar{E})$ has simply closed polygonal boundary,
\begin{equation} \label{GBsmooth}
   \sum_{x\in X} K(x) = 6 \chi(X).
\end{equation}
\end{theorem}
This result, and its proof, is a modest extension of recent work of Knill who establishes the discrete Gauss-Bonnet formula \eqref{GBsmooth} for subdomains of the triangular lattice (see \cite{Knill1} Section 7, and see \cite{Knill2} for generalizations to higher dimension). The extra term appearing in \eqref{GBrough} for irregular graphs, which we have not seen previously in the literature, may be viewed as a topological characteristic of ``boundary defects'' such as the wire edges depicted in Figure \ref{F:edge}.  
\\[2mm]
{\bf Proof} We use the second of the following three ``handshake'' properties (terminology in \cite{KnillSlides}) for any planar triangular graph $(X,\bar{E})$. Letting
\begin{equation}\label{Puiseux0}
  V_{0}(x):= 1, \;\; V_1(x):=\sharp\mbox{(edges in S(x))}, \;\;
  V_2(x) := \sharp\mbox{(faces in S(x))}
\end{equation}
and 
$$
  v_0 := \sharp\mbox{vertices in $X$}, \;\; v_1 := \sharp\mbox{edges in $X$}, \;\;
  v_2 := \sharp \mbox{faces in $X$},
$$
we claim that
\begin{enumerate}
\item $\sum\limits_{x\in X} V_0(x) = v_0,$
\item $\sum\limits_{x\in X} V_1(x) = 2 v_1 - \sharp  \mbox{(regular bdry edges)} - 2 \sharp \mbox{(wire edges)} = 2\sharp \bar{E} -  P(X)$,
\item $\sum\limits_{x\in X} V_2(x) = 3 v_2$.
\end{enumerate}
The first property is trivial. The second one follows because each interior edge appears as an edge in two unit spheres, each regular boundary edge appears in one unit sphere, and each wire edge does not appear in any unit sphere. And the third property follows from the fact that each face is triangular, and hence appears in precisely three unit spheres.

The following Euler-like curvature will emerge naturally during the computation below. It is reminiscent of \eqref{Euler} but uses the quantities \eqref{Puiseux0} instead of the edge and face degree (and we have changed normalization to obtain an integer quantity): 
\begin{equation}\label{Euler-like}
  \tilde{K}(x) = 6 \Bigl(V_0(x)-\frac12 V_1(x) + \frac13 V_2(x)\Bigr) = 6 - \sharp\mbox{edges in }S(x).
\end{equation}
We compute
\begin{eqnarray*}
  \chi(X) &=& \sum_{k=0}^2 (-1)^k v_k \;\; \mbox{(definition of Euler characteristic)} \\
          &=& \sum_{k=0}^2 (-1)^k \sum_{x\in X} \frac{V_k(x)}{k+1} \; - \; \frac{1}2  P(X) 
              \;\; \mbox{(handshake)} \\
          &=& \sum_{x\in X} \sum_{k=0}^2 (-1)^k \frac{V_k(x)}{k+1} \; - \; \frac{1}2  P(X)         
              \;\; \mbox{(change order of summation)} \\
          &=& \sum_{x\in X} \frac{\tilde{K}(x)}{6} - \frac{1}2  P(X) 
              \;\; \mbox{(definition of Euler-like curvature)} \\
          &=& \sum_{x\in X} \frac{K(x)}{6} + \frac{1}2 \Bigl(\sharp \, \partial X -  P(X)\Bigr) 
              \;\; \mbox{(definition of Puiseux curvature)}. 
\end{eqnarray*}
This establishes the theorem. Note how the surface term in the last but one line is cancelled, at least in the case of graphs with simply closed polygonal boundary, by the passage from Euler-like curvature to Puiseux curvature. 
\section{Geometric energy decomposition} \label{S:energy}
We now partition the at first sight featureless energy \eqref{energy} into various geometric and topological contributions. The pair potential may be any function $V\, : \, (0,\infty)\to \R\cup\{+\infty\}$ which attains its minimum at $r=1$ and has minimum value $-1$. 
We begin with a trivial but useful decomposition into minus the number of bonds of the original bond graph $(X,E)$, elastic energy, and the energy due to non-bonded interactions. Letting
\begin{eqnarray} \label{Eel}
   & & \calE_{e\ell}(X)  :=  \frac12 \sum_{\alpha\le |x-y|\le \beta} \Bigl(V(|x-y|)-\min_{r>0} V(r)\Bigr), \\
   & & \calE_{nbond}(X)  :=  \frac12 \sum_{|x-y| \not\in [\alpha,\beta]} V(|x-y|),
                 \label{Enbond} 
\end{eqnarray}                 
we have
\begin{equation} \label{decomp0}
    \calE_{_V}(X) = - \sharp {E} + \calE_{e\ell}(X) + \calE_{nbond}(X).
\end{equation}
Note that the elastic energy $\calE_{e\ell}$ vanishes unless there are stretched or compressed bonds; it is caused by particle pairs which are bonded (i.e., $\{x,y\}$ belongs to the edge set $E$) but whose distance is not the optimal bond length $r=1$.

The first term in \eqref{decomp0}, i.e. the number of edges in the untriangulated bond graph, can be decomposed further. First we introduce the following {\it defect measure} which quantifies 
the distance of the bond graph $(X,E)$ from a triangular graph:
\begin{equation} \label{mu}
  \mu(X)  :=  \sharp \bar{E} - \sharp E \;\;\;(= \mbox{number~of~additional~edges~due~to~triangulation}) 
\end{equation}
so that 
\begin{equation} \label{decomp02}
    -\sharp E = - \sharp \bar{E} + \mu(X).
\end{equation}
The first term on the right hand side of \eqref{decomp02} has sufficient regularity to be amenable to differential-geometric analysis. Recall the definitions \eqref{Puiseux0} and the associated ``handshake'' properties from subsection \ref{S:DGB}. We calculate
\begin{eqnarray}
  - \sharp\bar{E} & = & - \frac{1}{2} \Bigl( \sum_{x\in X} V_1(x) + P(X)\Bigr) \;\; \mbox{(by handshake)} \nonumber \\
                  & = & - 3 \, \sharp\, \mbox{int}\, X - \frac{3}{2} \sharp \, \partial X
                        + \frac{1}{2}\Bigl(\sum_{x\in X} K(x) -  P(X) \Bigr) \;\; 
\mbox{(by definition of $K$)} \nonumber \\
                  & = & - 3 N + \Bigl(\frac{3}{2} \sharp\, \partial X - \frac{1}{2} P(X) \Bigr) 
                          + \frac{1}{2} \sum_{x\in X} K(x) \;\;
(\mbox{since }\sharp X = N). 
                          \label{decomp1}
\end{eqnarray}
Combining \eqref{decomp0}, \eqref{decomp02}, \eqref{decomp1}, and the discrete Gauss-Bonnet theorem gives the following final energy decomposition which we state as a theorem. Note that the explicit contribution from boundary {\it particles} in \eqref{decomp1} and that from the curvature sum exactly cancel, even for irregular configurations, so that the surface energy contribution reduces to the perimeter \eqref{peri}. 
\begin{theorem} \label{T:geo} (Geometric energy decomposition) Let 
$V\, : \, (0,\infty)\to \R\cup\{+\infty\}$ be any pair potential which attains its minimum at $r=1$ and has minimum value $-1$.
For any number $N$ of particles, any N-particle configuration $X$ satisfying the mild minimum distance bound \eqref{planar2}, and any choice of the bond range $[\alpha,\beta]$ satisfying \eqref{planar1}, the atomistic energy \eqref{energy} satisfies the following exact decomposition:
\begin{equation} \label{decomp2}
 \calE_V(X) = -3 
N 
+ P(X) + 3\chi(X) + \mu(X) + \calE_{e\ell}(X) + \calE_{nbond}(X).
\end{equation}
Here $P(X)$ is the perimeter \eqref{peri}, $\chi(X)$ the Euler characteristic \eqref{Euchar}, $\mu(X)$ the defect measure introduced in \eqref{mu}, $\calE_{e\ell}$ the elastic energy \eqref{Eel}, and $\calE_{nbond}$ the energy \eqref{Enbond} due to non-bonded interactions.                   
\end{theorem}
The different terms in this decomposition are indicative of the different mechanisms of energy lowering which lead to crystallized ground states of special shape: lowering the Euler characteristic $\chi(X)$ aggregates different connected components into a single one; lowering the perimeter $P(X)$ condenses long chains into bulk regions; lowering the defect measure $\mu$ means forming more and more triangular faces; lowering the elastic energy means rigidifying stretched or compressed bonds into unit-length ones. 

We now comment on some of the individual terms.

The two terms $\mu(X)$ and $\chi(X)$ were introduced with the help of triangulation, but they are in fact independent of the triangulation. Regarding $\chi$ this was explained in Section \ref{S:IandB}, and regarding $\mu$ this follows from the following expression in terms of the faces $f$ of the original bond graph $(X,E)$:
\begin{equation} \label{mudirect}
  \mu(X) = \sum_f \Bigl( P^{inn}(f) - 3\Bigr), 
\end{equation}
where $P^{inn}(f)$ is the ``inner perimeter'' of the face $f$ which, analogously to the outer perimeter \eqref{peri}, counts wire edges twice,
\begin{equation} \label{Pinn}
  P^{inn}(f) = \sharp \, \mbox{regular interior edges in }\partial f 
             + 2\,\sharp \, \mbox{inner wire edges in }\partial f.   
\end{equation}
In particular, for connected $X$, $\mu(X)$ is the sum of the inner perimeter of the non-triangular faces minus 3 times the first Betti number $b_1$ of the closed subset in the plane formed by the union of the vertices, edges, and triangular faces.  

Finally we emphasize the following point. For typical configurations below the crystallization temperature, the first two terms, combinatorial bulk energy and combinatorial surface energy, have a clear asymptotic scaling with respect to the particle number $N$, namely $\sim N$ and $\sim N^{1/2}$, but the defect measure $\mu$ and the elastic energy $\calE_{e\ell}$ do {\it not}. Hence the above decomposition cannot be discovered via asymptotic analysis. 
Instead, it can be used to study different regimes, such as a single dislocation ($\mu \sim 1$,  $\calE_{e\ell}\sim \log N$), a dilute dislocation density as relevant for plasticity \cite{GLP} ($\mu\sim N/(\log N)^2$), or lattice elasticity ($\mu=0$ and ${\calE}_{e\ell}\sim N$). 
\section{Adding and removing a closed shell} \label{S:shells}
\begin{figure}[http!]
\begin{center}
\includegraphics[width=0.4\textwidth]{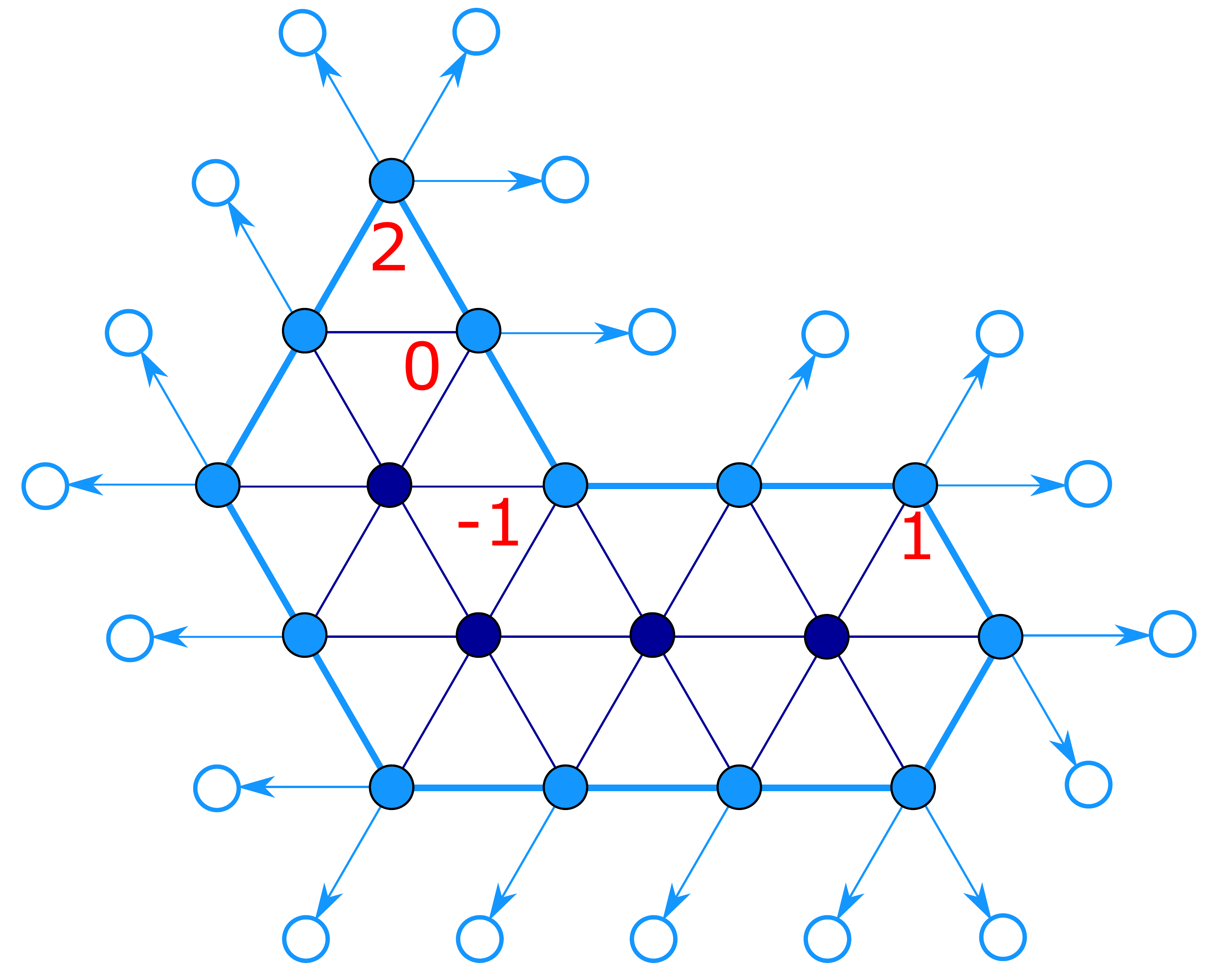}
\end{center}
\caption{Adding a closed shell to a crystallized configuration $X'$. The size of this closed shell can be counted using discrete curvature: for each particle $x\in\partial X'$ we need to add precisely $K(x)+1$ particles, where $K(x)$ is the combinatorial Puiseux curvature (see \eqref{Puiseux}; values shown in red).} 
\label{F:upper}
\end{figure}

\begin{figure}[http!]
\begin{center}
\includegraphics[width=0.6\textwidth]{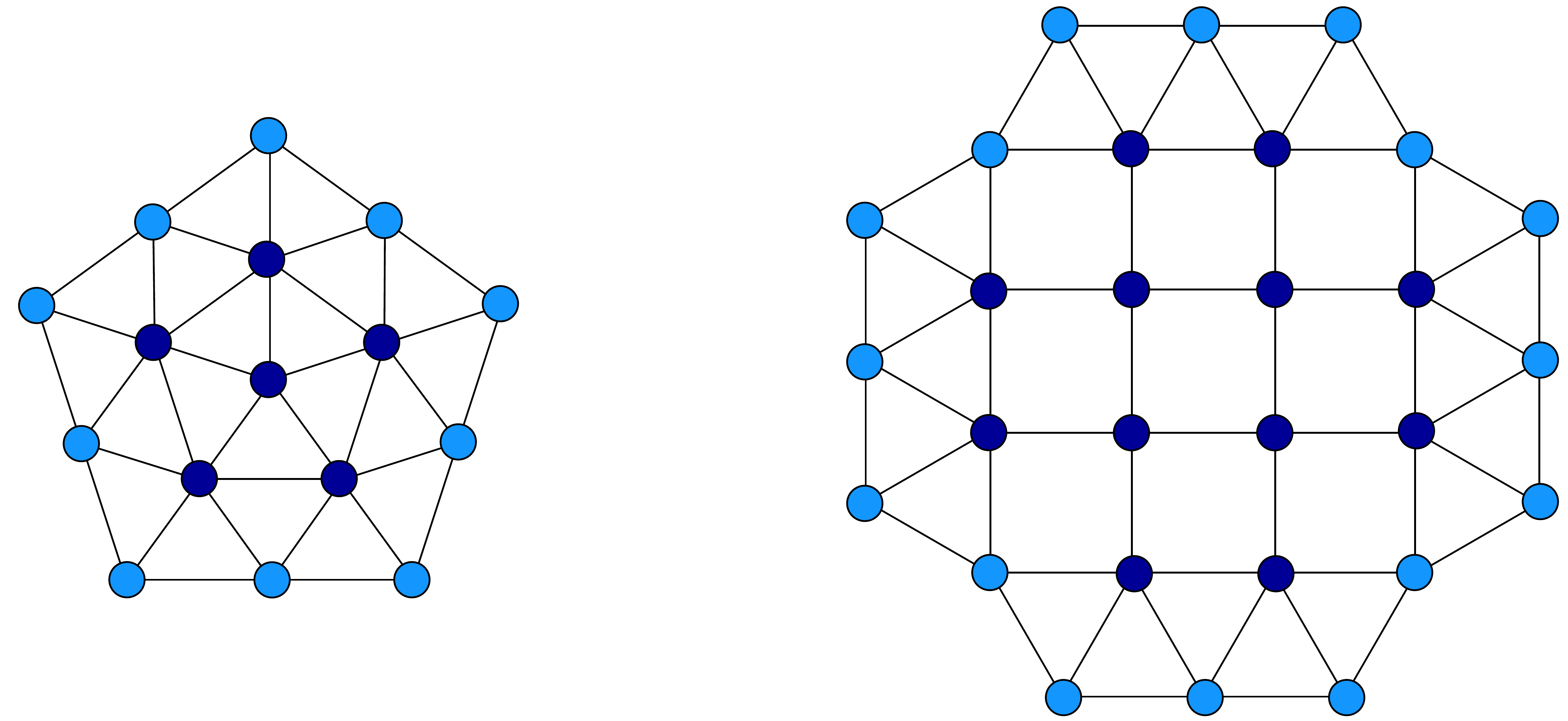}
\end{center}
\caption{A central phenomenon which makes crystallization hard to prove: removing the boundary from a noncrystallized configuration can reduce the perimeter by fewer than 6 edges. Left: an example with unequal bond length where $P(X)=10$, $P(X\backslash\partial X)=5$. Right: an example with equal bond lengths where $P(X)=16$, $P(X\backslash\partial X)=12$. Note that locally near the boundary, both examples look similar to crystallized configurations.} 
\label{F:counterex}
\end{figure}

An important idea of Heitmann and Radin \cite{HR} in their proof of crystallization for the potential \eqref{HR} (see also \cite{Harborth}) was to analyze the energy change when adding or removing not a single atom, but a {\it complete layer}, or ``closed shell'', of boundary atoms. This energy change can be derived in a transparent way using discrete curvature and Gauss-Bonnet, as we show in this section.
\\[2mm]
We begin by looking at crystallized configurations, i.e. subsets of the triangular lattice $\calL$. Throughout this section, in the definition of the bond graph of a configuration we require edges to be of euclidean length $1$, i.e. we choose $\alpha=\beta=1$ in \eqref{edge}.  
\begin{lemma} \label{L:closedshell} (Adding a closed shell to a crystallized configuration)
Let $X'$ be a crystallized configuration (i.e., a subset of the triangular lattice $\calL$)
with simply closed polygonal boundary. 
Let $X$ be the crystallized configuration obtained by adding every point in $\calL$ linked to $X'$ by an edge, i.e. 
\begin{equation} \label{cs}
   X := X' \cup \{x\in \calL\, : \, \mbox{there exists }x'\in X' \mbox{ with }|x-x'|=1\}. 
\end{equation}
Then the perimeter of $X$ satisfies the upper bound
\begin{equation} \label{upper}
   P(X) \le P(X') + 6.
\end{equation}
Moreover, provided, e.g., $X'$ contains at most two points with negative curvature, equality holds in \eqref{upper}, and the set $X$ also has a simply closed polygonal boundary and at most two points with negative curvature.   
\end{lemma}
{\bf Proof} Since $X'$ is crystallized with simply closed polygonal boundary, the points in $\partial X'$ have curvature $K=-1$, $0$, $1$, or $2$. See Figure \ref{F:upper}. 
Choose an orientation of the boundary polygon formed by joining the particles in $\partial X$. 
For any $x\in\partial X'$, we add $(K(x)+1)$ new particles with angles $j \pi/3$ 
($j=1,..,K(x)+1$) with respect to the incoming edge as shown in Figure \ref{F:upper}. 
This yields all particles in $X\backslash X'$, and moreover yields them exactly once provided $X'$ contains at most 2 points with negative curvature. It follows that
\begin{eqnarray*}
   \sharp(X\backslash X') & \le & \sum_{x\in\partial X'} (K(x)+1) \;\;(\mbox{with equality if }
                                  \sharp\{x\in X'\, : K(x)<0\}\le 2) \\
                          & = & \sum_{x\in X'} K(x) + \sharp \, \partial X' \;\;(\mbox{since $K(x)=0$ in int}\, X') \\
                          & = & 6 \chi(X') + \sharp \, \partial X' \;\;(\mbox{by Gauss-Bonnet}).                           
\end{eqnarray*}
Since $X'$ has simply closed polygonal boundary, we have $\chi(X')=1$ and $\sharp \, \partial X'=P(X')$. Moreover since the additional particles in $X$ have distance $1$ from $X'$, $X\backslash X'= \partial X$ and $X$ has no wire edges. In particular, $\sharp(X\backslash X') = \sharp \, \partial X = P(X)$. This establishes \eqref{upper} and the fact that equality holds when $X'$ has at most two points of negative curvature. The additional consequences of the latter property are straightforward. The proof of the lemma is complete.
\\[2mm]
For typical crystalline configurations, it is easy to convince oneself that the reverse operation -- removing the boundary -- lowers the perimeter by 6.
A new phenomenon occurs for non-crystalline configurations: {\it removing the boundary can reduce the perimeter by less than 6}. This phenomenon is in our view the central reason why the crystallization problem in two dimensions is difficult even for the HR potential \eqref{HR}, and why it remains open for short-range soft potentials. For examples see Figure \ref{F:counterex}.
\vspace*{2mm}

\noindent
The key to overcome this phenomenon in case of the potential \eqref{HR} is the lemma below, which restores the number 6 by accounting for defects. 
\begin{lemma} \label{L:lower} (Removing the boundary from a hard-sphere configuration) 
Let $X$ be a hard-sphere configuration, that is to say a configuration in $\R^2$ with the two properties that \\
(i) all euclidean interparticle distances are $\ge$ 1, \\
(ii) all edges have euclidean length 1 (i.e., $\alpha=\beta=1$ in \eqref{edge}). \\
Assume also that $X$ has simply closed polygonal boundary, and that $X'=X\backslash\partial X$ is nonempty. Then 
\begin{equation} \label{lower}
   P(X) + \mu(X) \ge P(X') + \mu(X') + 6.
\end{equation}
Moreover if equality holds, then we must have $\mu(X)=\mu(X')$. 
\end{lemma}
The configuration in the right part of Figure \ref{F:counterex} shows that \eqref{lower} can fail if the defects terms $\mu(X)$ and $\mu(X')$ are dropped.  
\\[2mm]
The key to the proof is the following curvature bound for hard-sphere configurations. 
To state it we first need to deal with a certain ``arbitrariness'' that entered by triangulation. The triangulation will typically introduce extra edges emanating from some of the boundary particles, and therefore influence the Puiseux curvature. However, for configurations with simply closed polygonal boundary, given any fixed $x\in \partial X$ there always exists a triangulation $(X,\bar{E})$ such that no additional edges emanate from $x$. This suggests to consider the following triangulation-independent variant of the Puiseux curvature \eqref{Puiseux} 
\begin{eqnarray} \label{Kmax}
          K_{max}(x) & := & \max \{ K(x)\, : \, (X,\bar{E}) \mbox{ is a triangulation of }(X,E)\} \\
                     & =  & 2 - i(x),       \nonumber      
\end{eqnarray}
where here and below $i(x)$ denotes the number of interior edges in the original bond graph $(X,E)$ which emanate from $x$. 
\\[2mm]
\begin{proposition} \label{P:curvaturebound} (Discrete curvature bounds euclidean curvature) \\
Let $X$ be a two-dimensional hard-sphere configuration (see Lemma \ref{L:lower}) with simply closed boundary. Then the maximal Puiseux curvature \eqref{Kmax} of the boundary (with respect to any triangulation $(X,\bar{E})$ of the bond graph) is bounded from below by the euclidean Puiseux curvature \eqref{Keu}: 
\begin{equation}\label{curvaturebound}
                  \frac{K_{max}(x)}{|S_{1,\calL}|} \ge \frac{K_{eu}(x)}{|S_1|_{eu}} \;\;\mbox{ for all }x\in\partial X.
\end{equation}
Here $|S_{1,\calL}|=6$ is the length w.r. to the graph metric of the unit sphere in the triangular lattice, and $|S_1|_{eu}=2\pi$ is the euclidean length of the euclidean unit sphere in $\R^2$. 
\end{proposition}
{\bf Proof of Proposition \ref{P:curvaturebound}} Let $x$ be any boundary particle in $X$. Because of the euclidean distance constraint on the particles in $X$, if $i(x)$ interior edges emanate from $x$, then the inner angle $\alpha(x)$ between the incoming and outgoing boundary edge at $x$ (see Figure \ref{F:Puiseux}) must satisfy 
\begin{equation} \label{anglebd1}
   \alpha(x) \ge \Bigl( i(x) + 1\Bigr) \cdot \frac{\pi}{3},
\end{equation}
as already noted by Harborth \cite{Harborth}. We obtain the differential-geometric meaning of 
\eqref{anglebd1} by subtracting $\pi$ from both sides and multiplying by $-\frac{1}{2\pi}$: 
\begin{equation} \label{anglebd2}
   \frac{K_{eu}(x)}{|S_1|_{eu}} = \frac{\pi - \alpha}{2\pi} \le \frac{2- i(x)}{6}
   = \frac{K_{max}(x)}{|S_{1,\calL}|}.  
\end{equation}
\\[2mm]
{\bf Proof of Lemma \ref{L:lower}} Introduce the subset of interior edges in the (non-triangulated) bond graph which touch the boundary,
$$
   I := \{\mbox{interior edges of }(X,E) \, \mbox{with at least one vertex in }\partial X\}. 
$$
Moreover, abbreviate the Heitmann-Radin energy $\calE_{\VHR}$ defined by 
\eqref{energy}--\eqref{HR} by $\calE_{HR}$. The energy difference $\calE_{HR}(X) - \calE_{HR}(X')$ can be computed in 2 ways. On the one hand, it trivially equals the number of of edges in the bond graph $(X,E)$ with at least one vertex in $\partial X$, i.e.
$$
   \calE_{HR}(X) - \calE_{HR}(X') = - P(X) - \sharp I.
$$
On the other hand, applying Theorem \ref{T:geo} to both $X$ and $X'$ gives
$$
  \calE_{HR}(X) - \calE_{HR}(X') = - 3 \sharp \, \partial X + (P+\mu)(X) - (P+\mu)(X') + 3(\chi(X)-\chi(X')). 
$$
Equating both expressions yields 
\begin{equation} \label{magic}
    (P+\mu )(X) - (P+\mu )(X') = 2 \, \sharp \, \partial X - \sharp I + 3\underbrace{(\chi(X')-\chi(X))}_{\ge 0}. 
\end{equation}
We claim that the underbraced term is nonnegative, i.e. that $\chi(X')\ge 1$.
This can be seen as follows. Because the original graph was assumed connected, the union of the closure of all its faces equals the complement of the ``outer'' region (i.e., the unbounded connected component of the complement). Hence the connected components of $X'$ are all simply connected, so the Euler characteristic of $X'$ equals the no. of connected components of $X'$.

To maximise clarity of the argument below, we split the set $I$ of interior vertices touching $\partial X$ into the two sets
\begin{eqnarray*}
   & & I_1 := \{ e\in I \, : \, \mbox{precisely 1 vertex of $e$ belong to }\partial X\}, \\
   & & I_2 := \{ e\in I \, : \, \mbox{both vertices of $e$ belong to }\partial X\}.
\end{eqnarray*}
It follows that 
\begin{eqnarray} 
   \frac{(P+\mu )(X)-(P+\mu)(X')}{|S_{1,\calL}|} & \ge &  \frac{2\, \sharp \, \partial X - \sharp I}{|S_{1,\calL}|} \;\; \mbox{(by \eqref{magic})} \nonumber \\
            & \ge &  \frac{2\, \sharp\,\partial X - (\sharp I_1 + 2 \sharp I_2)}{|S_{1,\calL}|} \;\; \mbox{(trivially)} \nonumber \\
            & = &  \frac{\sum_{x\in\partial X} K_{max}(x)}{|S_{1,\calL}|} \;\; \mbox{(by \eqref{Kmax})} \nonumber \\
            & \ge &  \frac{\sum_{x\in\partial X} K_{eu}(x)}{|S_1|_{eu}} \;\; \mbox{(by the curvature bound \eqref{curvaturebound})} \nonumber \\
            & = & 1 \;\; \mbox{(by classical Gauss-Bonnet)}. \label{chain}
\end{eqnarray}

In the last line we have used the classical Gauss-Bonnet theorem in the plane which says that the curvature $\kappa$ of a simply-closed smooth planar curve $\calC$ satifies 
\begin{equation} \label{GBclass}
   \int_{\calC} \kappa \, ds = 2\pi,
\end{equation}
where $ds$ denotes the line element (alias Hausdorff measure $H^1$) on $\calC$. (The curvature $\kappa$ in \eqref{GBclass} is defined by the relation $\kappa(c(t))=\alpha'(t)$, where $c \, : \, [0,T]\to\R^2$ is a smooth anticlockwise parametrization of $\calC$ and $\alpha(t)$ is the angle of the tangent vector $c'(t)$ with respect to a reference direction, measured anticlockwise.) 
More precisely, we have used the limiting case where $\calC$ is the piecewise linear curve given by the union of the boundary edges of $X$ and the measure $\kappa\, ds$ on $\calC$ is the singular measure $\sum_{x\in\partial X} K_{eu}(x) \delta_x$. This establishes \eqref{lower}. 

When equality holds in \eqref{lower}, and hence in the last inequality in \eqref{chain}, then for any $x\in\partial X$ we must have equality in \eqref{anglebd1}; but this implies that all angles between any two consecutive edges emanating from $x$ except the outer angle between the two boundary edges are equal to $\pi/3$, and therefore the corresponding endpoints must be joined by an edge. This excludes the possibility that any additional edge in the triangulated bound graph can emanate from $x$. The proof of Lemma \ref{L:lower} is complete.   
\section{A new proof of crystallization in the Heitmann-Radin model}
\begin{figure}[http!]
\begin{center}
\includegraphics[width=0.8\textwidth]{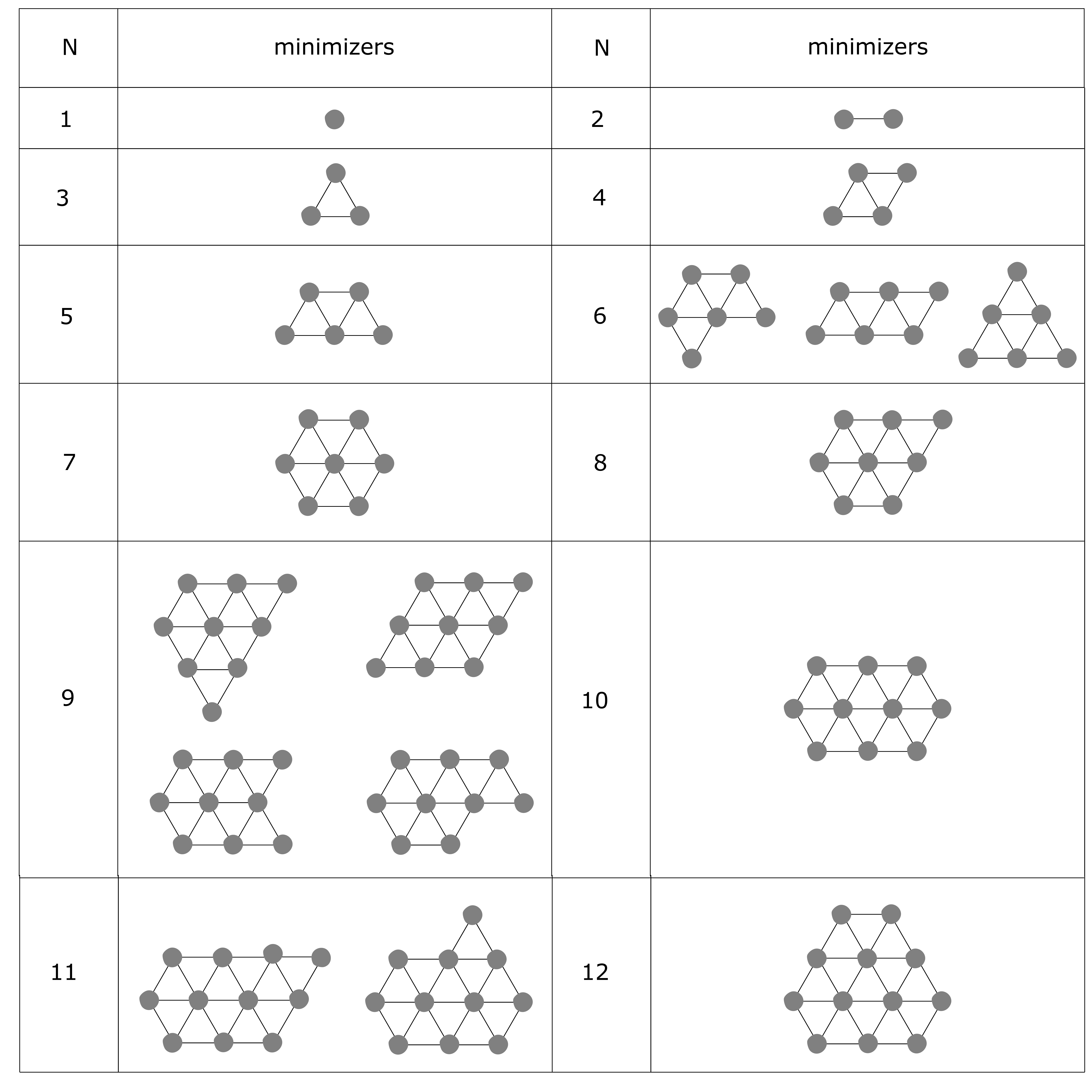}
\end{center}
\caption{The minimizers of the Heitmann-Radin energy up to $N=12$.} 
\label{table:minimizers}
\end{figure}

The preceding tools allow a novel proof of the Heitmann-Radin crystallization theorem. One observes that geometric rigidity of minimizers follows immediately from ``topological crystallization'' of the bond graph (defined below), and establishes the latter by control of the defect measure $\mu(X)$.
%
%
%
%
%
\begin{theorem} \label{T:HR} For any number $N$ of particles, minimizing configurations $X=\{x_1,..,x_N\}$ of the Heitmann-Radin energy \eqref{energy}, \eqref{HR} are, up to an overall rotation and translation, subsets of the triangular lattice \eqref{lattice}. For $N\ge 3$, minimizers have simply closed boundary. Moreover there always exists a minimizer containing at most one point with negative curvature.  
\end{theorem}
\noindent
{\bf Definition} {\it A configuration $X$ with bond graph $(X,E)$ is called topologically crystallized if it satisfies $\mu(X)=0$ and has simply closed boundary.}
\\[2mm]
This property implies the asserted geometric rigidity because -- as explained in Section \ref{S:shells} -- in the present context edges are of length $1$, whence all faces of a configuration with $\mu(X)=0$ are equilateral triangles.  
\\[2mm]
{\bf Proof} First we note that for $N\le 9$ the result holds, by inspection. See Figure \ref{table:minimizers} for the minimizers in these cases. We will proceed by induction and assume that the result holds for $N$ up to some particle number $N_0\ge 9$. We now consider $N=N_0+1$, and suppose that $X$ is an $N$-particle minimizer. Note that $X$ must be connected and have simply closed boundary, because otherwise we could lower the energy by translating or rotation parts of the configuration against each other. This implies, in particular, that $P(X)=\sharp\, \partial X$. We need to show that $\mu(X)=0$. 

{\bf Step 1: Competitor.} We construct a competitor of $Y$ of $X$ with the same number of particles, as follows:
\begin{equation} \label{competitor}
  X \xrightarrow{\mbox{\rm\scriptsize remove boundary}} X' \xrightarrow{\mbox{\rm\scriptsize ground state with }\sharp X' \mbox{\rm\scriptsize  particles}} Y' \xrightarrow{\mbox{\rm\scriptsize add $\sharp\,\partial X$ particles}} Y.
\end{equation}
The first two maps are easily made precise: $X' := X \backslash \partial X$, and $Y'$ is any minimizer subject to the number of particles being $\sharp X'$. By the induction hypothesis, $Y'$ is a subset of the triangular lattice, has simply closed boundary, and can be chosen to contain at most one point with negative curvature. To describe the third map in \eqref{competitor} is a little more work, since we need to add precisely $d:=\sharp \, \partial X$ particles to $Y'$, to restore the original particle number $N$. We do so using the following unique decomposition of any such $d\ge 0$ into the number of particles in as many additional closed shells (see \eqref{cs}) as possible and a remainder,
\begin{equation} \label{remainder}
  d = d_m + \delta < d_{m+1}, \;\;\; m\in\{0,1,2,...\}, \;\;\; \delta \in\{0,1,2,...\}.
\end{equation}
Hence to obtain the competitor $Y$ in \eqref{competitor} we add $m$ closed shells to $Y'$, resulting in a configuration $Y_m$, and then add $\delta$ additional particles from the next shell $Y_{m+1}\backslash Y_m$ which form a connected subset thereof. (If $\sharp X'=0$, we set $Y'=\emptyset$ and let the first closed shell around it be the first 6-particle minimizer in Table \ref{table:minimizers}.) It is straightforward to verify using induction over $m$ and Lemma \ref{L:closedshell} that
\begin{equation}
   P(Y_m) = P(Y') + 6m, \;\;\; d_m(Y') = m \Bigl(P(Y')+3(m+1)\Bigr).
\end{equation}
The remaining $\delta$ particles are added
starting with a neighbour of a point with negative curvature if such a point exists, or otherwise at the end of the longest straight side. More precisely, in the latter case we start with a joint neighbour of the first two particles belonging to this side. We note that adding $\delta$ additional particles in this way increases the perimeter at most by $\delta$, and at most by $\delta-1$ when $\delta\ge 2$. The improved bound for $\delta\ge 2$ relies on the fact that $N\ge 10$, because, e.g., for $N=9$, the anomalous situation can occur that $Y'$ consists of a single point and $Y_1$ of the hexagon-shaped 7-particle minimizer depicted in Table \eqref{table:minimizers}, in which case the perimeter increases by $\delta$ when $\delta=2$. 

{\bf Step 2: Absence of defects.} The energy decomposition \eqref{decomp2} now provides a powerful tool, because we can separately compare the individual energy contributions. By \eqref{intro:decomp1} and the fact that $X$ and $Y$ have the same number of particles, we only need to compare Euler characteristic, defect measure, and perimeter. The first two a trivial to deal with: 
\begin{equation} \label{muandchi}
   0=\mu(Y)\le\mu(X), \;\;\; 1=\chi(Y) = \chi(X).
\end{equation}
To compare the perimeter is less trivial. We distinguish three cases, corresponding to $X$ having a large, moderate, or small number of boundary particles, where large/small means large/small compared to the number $d_1(Y')$ of particles of a single closed shell around $Y'$. 
\\[1mm]
a) {\it Large number of boundary particles:} $\sharp \, \partial X \ge d_1(Y') + 2$. \\
This means that either $m\ge 2$ or $m=1$, $\delta\ge 2$. In the first case, $Y$ is obtained by adding $m\ge 2$ closed shells and $\delta\ge 0$ additional particles. By the fact that the additional particles increases the perimeter at most by $\delta$, Lemma \ref{L:closedshell}, and the elementary to check inequality that $P(Y')+6m < d_m$ when $m\ge 2$, we have 
$$
   P(Y) \le P(Y_m) + \delta = P(Y') + 6m + \delta < d_m(Y') + \delta = \sharp \, \partial X = P(X). 
$$
In the second case, $Y$ is obtained by adding one closed shell and $\delta\ge 2$ additional particles. In this case $P(Y')+6m = d_1$ but, as explained at the end of Step 2, the $\delta$ additional particles raise the perimeter by at most $\delta-1$, so we obtain the same overall inequality:
$$
  P(Y) < P(Y_m) + \delta = P(Y') + 6 + \delta = d_1(Y') + \delta = \sharp \, \partial X = P(X). 
$$
Hence in both cases we have $P(Y) < P(X)$, which 
together with \eqref{muandchi} shows that $X$ is not a minimizer, contradicting our assumption. Hence case a) cannot occur. 
\\[1mm]
b) {\it Moderate number of boundary particles:} $\sharp \, \partial X = d_1(Y') + 1$. \\
In this case $Y$ is obtained from $Y'$ by adding one closed shell and attaching one additional particle. Using the trivial fact that the additional particle raises the perimeter at most by 1, Lemma \ref{L:closedshell}, and the formula for $d_1$ we estimate
\begin{equation} \label{perestb}
  P(Y) \le P(Y_1) + 1 = P(Y')+6 + 1 = d_1(Y') + 1 = \sharp \, \partial X = P(X).
\end{equation}
Since $X$ is a minimizer, this together with \eqref{muandchi} 
shows that $\mu(X)=0$.
\\[1mm]
c) {\it Small number of boundary particles:} $\sharp \, \partial X = d_1(Y')$. \\
In this case $Y$ is obtained from $Y'$ by adding one closed shell, or a 
connected part thereof. Using, in order of appearance, that adding only a connected part thereof does not give higher perimeter than adding the full closed shell, Lemma \ref{L:closedshell}, the fact that $Y'$ is defect-free by the induction hypothesis, the minimizing property of $Y'$, and Lemma \ref{L:lower} yields
\begin{multline}\label{perestc}
  P(Y) \le P(Y_1) = P(Y')+6 = P(Y')+ \mu(Y') +6 \le (P+\mu)(X') + 6 \\
\le (P+\mu)(X).
\end{multline}
Since $X$ is a minimizer, this together with \eqref{muandchi} shows that all inequalities above are equalities. But by Lemma \ref{L:lower}, equality in the last inequality implies $\mu(X)=\mu(X')$. And equality in the last but one inequality implies that $X'$ is a minimizer, whence $\mu(X')=0$, by the induction hypothesis. Combining these two statements gives $\mu(X)=0$. 

Moreover we have shown in both case b) and case c) that the competitor $Y$ is also a minimizer. Since by construction $Y$ has at most one point with negative curvature, this establishes the ``there exists'' statement in the theorem, completing the proof. 
\section{Concluding remarks}
Perhaps the main advance in our derivation of Heitmann-Radin crystallization is that the exact numerical value of the ground state energy is no longer needed. Instead, proving crystallization is decoupled from minimizing the ensuing lattice model; note that for lattice models, powerful tools are available, including at finite temperature (see e.g. \cite{DKS}). 
\\[4mm]               
{\bf Acknowledgements.}
This work was supported by the DFG Collaborative Research Center TRR 109 ``Discretization in Geometry and Dynamics''. GF thanks Sasha Bobenko and John Sullivan for their most helpful advice at an early stage of this project that discrete Gauss-Bonnet theorems might be relevant to the goal of relating defect-induced curvature to atomistic energies. Also, LDL and GF thank Oliver Knill for sharing valuable intuition and insights related to Refs. \cite{Knill1, Knill2}.


\begin{thebibliography}{99}
\bibitem{AFS} Au Yeung Y., Friesecke G., Schmidt B.: Minimizing atomic configurations of short range pair potentials in two dimensions: crystallization in the Wulff shape, \emph{Calc. Var. Partial Differential Equations}, \textbf{44} (2012), 81--100.



\bibitem{AllenTildesley} M.P. Allen, D.J. Tildesley: Computer Simulation of Liquids. \emph{Oxford University Press}, Oxford (1987).

\bibitem{CS} J.H. Conway, N.J.H. Sloane, \emph{Sphere Packings, Lattices and Groups} 3rd edn., Springer-Verlag, Berlin (1999).

\bibitem{CameronVandenEijnden} M. Cameron, E. Vanden-Eijnden, Flows in Complex Networks: Theory, Algorithms, and Application to Lennard-Jones Cluster Rearrangement, \emph{J. Stat. Phys.} 156, 3 (2014), 427--454.

\bibitem{DPS} E. Davoli, P. Piovano, U. Stefanelli, Sharp $N^{3/4}$ law for the minimizers of the edge-isoperimetric problem on the triangular lattice, preprint,  http://cvgmt.sns.it/media/doc/paper/2862/edgeperimeter13December2015.pdf.

\bibitem{DF2} L. De Luca, G. Friesecke, In preparation.

\bibitem{DKS}  Dobrushin R.L., Kotecky R., Shlosman S.B.: \emph{The Wulff Construction: A Global Shape from Local Interactions}, AMS: Providence (1992). 

\bibitem{ELi} W. E, D. Li, On the crystallization of 2D hexagonal lattices, \emph{Comm. Math. Phys.}, \textbf{286} (2009), 1099-1140.

\bibitem{Federer} H. Federer, Geometric measure theory

\bibitem{FTTT} Flatley L., Tarasov A., Taylor M., Theil F.: Packing twelve spherical caps to maximize tangencies, \emph{J. Comput. Appl. Math.}, \textbf{254} (2013), 220--225. 

\bibitem{FlatleyTheil} L.C. Flatley, F. Theil, Face-Centered Cubic Crystallization of Atomistic Configurations, {\it . Arch. Ration. Mech. Anal.}, \textbf{218} (2015), no. 1, 363--416.


\bibitem{FrieseckeTheil} G. Friesecke, F. Theil, Molecular Geometry Optimization, Models,
Encyclopedia of Applied and Computational Mathematics, 10.1007/978-3-540-70529-1$_-$239, Springer-Verlag (2015). 

\bibitem{GLP} A. Garroni, G. Leoni, M. Ponsiglione, Gradient theory for plasticity via homogenization of discrete dislocations, {\it J. Eur. Math. Soc.} {\textbf 12} (5) (2010), 1231--1266.

\bibitem{Gromov} M. Gromov, {\it Hyperbolic groups}, Essays in group theory, S. M. Gersten, (Editor), M.S.R.I. Publ. 8, Springer, 1987, 75--263.

\bibitem{Harborth} H. Harborth, L\"osung zu Problem 664A, \emph{Elem. Math.}, \textbf{29} (1974), 14--15.


\bibitem{HR} R.C. Heitmann, C. Radin, The ground states for sticky discs, \emph{J Stat. Phys.}, \textbf{22} (1980), no. 3, 281--287. 

\bibitem{Higuchi} Y. Higuchi, Combinatorial curvature for planar graphs, {\it Journal of Graph Theory}, \textbf{38} (2001), no.4,  220--229.

\bibitem{Knill1} O. Knill, A discrete Gauss-Bonnet type theorem, \emph{Elem. Math.}, \textbf{67}, 2012, 1--17.

\bibitem{Knill2} O. Knill, A graph theoretical Gauss-Bonnet-Chern theorem, \emph{arXiv} 1111.5395, 2011 

\bibitem{KnillSlides} O. Knill, Slides of talk at Joint Mathematics Meetings (JMM), Baltimore, 2014, available at
http://www.math.harvard.edu/~knill/seminars/baltimore/baltimore.pdf.



\bibitem{Radin} C. Radin, The ground states for soft discs, \emph{J Stat. Phys.}, \textbf{26} (1981), no.2, 281--287.



\bibitem{Schmidt} B. Schmidt, Ground states of the 2D sticky disc model: fine properties and $N^{3/4}$ law for the deviation from the asymptotic Wulff shape, \emph{J. Stat. Phys.}, \textbf{153} (2013),  727--738.

\bibitem{MaininiStefanelli} E. Mainini, U. Stefanelli, Crystallization in carbon nanostructures,
{\it Comm. Math. Phys.}, 328 (2014), no.2, 545--571.

\bibitem{Theil} F. Theil, A proof of crystallization in two dimensions, \emph{Comm. Math. Phys.}, \textbf{262} (2006), no. 1, 209--236.

\bibitem{TroianEtAl} P.A. Thompson, S.M. Troian, A general boundary condition for liquid flow ad solid surfaces, \emph{Nature} 389 (1997), 360--362. 

\bibitem{WalesDoye} D.J. Wales, Global optimization by basin-hopping and the lowest energy structures of Lennard-Jones clusters containing up to 110 Atoms, \emph{J. Phys. Chem. A}, 101 (1997), 5111--5116.

\end{thebibliography}
\end{document}